\documentclass[12pt]{amsart}

\usepackage[OT2, T1]{fontenc}

\usepackage{amscd}
\usepackage{amsmath}
\usepackage{amssymb}
\usepackage{mathrsfs} 			
\usepackage{units}
\usepackage[all]{xy}

\usepackage{algorithm}
\usepackage{algorithmic}

\def\frk{\frak}               

\def\mm{{\frk m}}

\def\Phi{{\frk n}}
\def\Phi{{\frk N}}
\def\opn#1#2{\def#1{\operatorname{#2}}} 
\opn\chara{char} \opn\length{\ell} \opn\pd{pd} \opn\rk{rk}
\opn\projdim{proj\,dim} \opn\injdim{inj\,dim} \opn\rank{rank}
\opn\depth{depth} \opn\sdepth{sdepth} \opn\fdepth{fdepth}
\opn\grade{grade} \opn\height{height} \opn\embdim{emb\,dim}
\opn\codim{codim}  \opn\min{min} \opn\max{max}

\opn\Tr{Tr} \opn\bigrank{big\,rank}
\opn\superheight{superheight}\opn\lcm{lcm}
\opn\trdeg{tr\,deg}
\opn\reg{reg} \opn\lreg{lreg} \opn\ini{in} \opn\lpd{lpd}
\opn\size{size}
\opn\div{div} \opn\Div{Div} \opn\cl{cl} \opn\Cl{Cl}
\opn\Spec{Spec} \opn\Supp{Supp} \opn\supp{supp} \opn\Sing{Sing}
\opn\Ass{Ass} \opn\Min{Min}
\opn\Ann{Ann} \opn\Rad{Rad} \opn\Soc{Soc}
\opn\Im{Im} \opn\Ker{Ker} \opn\Coker{Coker} \opn\Am{Am}
\opn\Hom{Hom} \opn\Tor{Tor} \opn\Ext{Ext} \opn\End{End}
\opn\Aut{Aut} \opn\id{id}  \opn\deg{deg}

\opn\nat{nat}
\opn\pff{pf}
\opn\Pf{Pf} \opn\GL{GL} \opn\SL{SL} \opn\mod{mod} \opn\ord{ord}
\opn\Gin{Gin} \opn\Hilb{Hilb}
\opn\aff{aff} \opn\con{conv} \opn\relint{relint} \opn\st{st}
\opn\lk{lk} \opn\cn{cn} \opn\core{core} \opn\vol{vol}
\opn\link{link} \opn\star{star}
\opn\gr{gr}

\def\pot#1#2{#1[\kern-0.28ex[#2]\kern-0.28ex]}

\opn\dirlim{\underrightarrow{\lim}}
\opn\inivlim{\underleftarrow{\lim}}
\let\to=\rightarrow

\def\Implies{\ifmmode\Longrightarrow \else
        \unskip${}\Longrightarrow{}$\ignorespaces\fi}
\def\implies{\ifmmode\Rightarrow \else
        \unskip${}\Rightarrow{}$\ignorespaces\fi}
\def\iff{\ifmmode\Longleftrightarrow \else
        \unskip${}\Longleftrightarrow{}$\ignorespaces\fi}

\let\:=\colon
\newtheorem{Theorem}{Theorem}[]
\newtheorem{Lemma}[Theorem]{Lemma}
\newtheorem{Corollary}[Theorem]{Corollary}

\theoremstyle{definition}

\newtheorem{Remark}[Theorem]{Remark}

\newtheoremstyle{subsection-tweak}
   {11pt}
   {3pt}%
   {}
   {}%
   {\bfseries}
   {}%
   {.5em}
   {\thmnumber{\@{#1}{}\@{#2}.}%
    \thmnote{~{\bfseries#3.}}}    

\newcounter{numberingbase}

\theoremstyle{subsection-tweak}
\newtheorem{bpp}[Theorem]{}
\newtheorem{bppt}[numberingbase]{}
\newcommand{\bbpp}{\begin{bpp}}
\newcommand{\eepp}{\end{bpp}}
\newcommand{\bbppt}{\begin{bppt}}
\newcommand{\eeppt}{\end{bppt}}

\theoremstyle{theorem}

\theoremstyle{definition}

\newcommand{\val}{\mathrm{val}}		

\let\epsilon\varepsilon
\let\phi=\varphi
\def\qed{\ifhmode\textqed\fi
      \ifmmode\ifinner\quad\qedsymbol\else\dispqed\fi\fi}
\def\textqed{\unskip\nobreak\penalty50
       \hskip2em\hbox{}\nobreak\hfil\qedsymbol
       \parfillskip=0pt \finalhyphendemerits=0}
\def\dispqed{\rlap{\qquad\qedsymbol}}

\opn\dis{dis}
\def\pnt{{\raise0.5mm\hbox{\large\bf.}}}

\opn\Lex{Lex}

\begin{document}

\title{A special form of Zariski's Uniformization Theorem in positive characteristic.}

\author{ Dorin Popescu}

\address{Simion Stoilow Institute of Mathematics of the Romanian Academy,
Research unit 5, P.O. Box 1-764, Bucharest 014700, Romania,}

\address{University of Bucharest, Faculty of Mathematics and Computer Science
Str. Academiei 14, Bucharest 1, RO-010014, Romania,}

\address{ Email: {\sf dorin.m.popescu@gmail.com}}

\dedicatory{To the $80^{th}$ birthday   of Vasile Br\^ inz\u anescu.}

\begin{abstract} This is mainly a small exposition on extensions of valuation rings
as a filtered union of smooth algebras. 
\end{abstract}

\maketitle
\section{Introduction}

A local ring   $(A,\mm)$ is  Noetherian if its  ideals are finitely generated. The maximum of the length of the chains
$p_0
\subsetneq p_1\subsetneq\ldots \subsetneq p_d$
 of its prime ideals is the Krull {\em dimension} of $A$. The Krull dimension of $A$ corresponds to dimension of algebraic varieties in algebraic geometry.
Another invariant is the  {\em embedding dimension} of $A$ which is 
 $\dim_{A/\mm} \mm/\mm^2$  as a linear space.

Let $C$ be the curve in the real plane
  ${\bf R}^2$ given by the equation   $X^2-Y^3=0$. The origin of the plane  $(0,0)$ is a point of  $C$, which is not smooth  in the usual terminology.  The dimension of the curve $C$ is $1$ in geometry  and its associated local ring  $B=({\bf R}[X,Y]/(Y^3-X^2))_{(X,Y)}$ in its point $(0,0)$ has a unique  big chain of prime ideals   $(0)\subset (X,Y)B$, that is $\dim B=1$. Note that 
  $\dim_{\bf R} (X,Y)B/(X,Y)^2B=2$.
on the other hand, the local ring   
   $B'=({\bf R}[X,Y]/(Y^3-X^2))_{(X-1,Y-1)}$ associated to  $C$ in its point $(1,1)$ has also a unique  big chain $(0)\subset (X-1,Y-1)B'$, that is the Krull dimension of   $B'$ is still $1$. Applying  Taylor's formula  for  $f=X^2-Y^3$ we get 
$$f=2(X-1)+(X-1)^2-3(Y-1) -3(Y-1)^2-(Y-1)^3$$ and so in the linear space 
  $( X-1,Y-1)B'/( X-1,Y-1)^2B'$ over  $\bf R$ there exists a linear dependence 
 between the elements induced by  
    $X-1,Y-1$ and so  the embedding dimension of  $B'$ is $1$.
 In general we have 
 $\dim A\leq \dim_{A/\mm} \mm/\mm^2 $ for a Noetherian local ring $(A,\mm)$. 
Note that this inequality is strict for $B$ but equality for $B'$.

 A Noetherian local ring $(A,\mm)$ is {\em regular} if the above inequality is an equality.
In geometry the local regular rings correspond to the smooth points of a variety. This is a reason that  $(1,1)$ is a smooth point of $C$ but $(0,0)$ not.

Consider the transformation 
   ${\bf R}\to C$ given by  $u\to (u^3,u^2)$ and we see that so we obtain a desingularization $\bf R$ of $C$, the real line being smooth.
  Algebraically this means to embed the domain  
    $A={\bf R}[X,Y]/(X^2-Y^3)$ in the subring $B$ generated by   $X/Y$ in the fraction field    $K$ of  $A$, which is finite over $A$ because     $(X/Y)^2-Y$ is $0$ in $K$. Note that $A$ and $B$ have the same fraction field. Actually, $B$ is the integral closure of $A$ in $K$ and we may use the integral closure to define the desingularization of curves even when we replace $\bf R$  
 by a field of positive characteristic.
 
    The desingularization of an algebraic variety of higher dimension was a hard problem of the last century solved by H.Hironaka in 1964 in the case of characteristic zero. His proof followed the 
Zariski Uniformization theorem \cite{Z}, which says that 
     a valuation ring V containing a field $K$ of characteristic zero is a filtered union of its regular  $K$-subalgebras. In positive characteristic such theorem is missing.
     
     An  $A$-algebra $D$ is  {\em
smooth} if $D$ is a localization of an $A$-algebra of type  $(A[Y]/(f))_M$, where  $f=(f_1,\ldots ,f_r) $  is a system of polynomials
in  $Y=(Y_1,\ldots,Y_N)$, $N\geq r$ over  $A$ and  $M$ is an $r\times r$-minor of $(\partial f/\partial Y)$. If $A$ is a regular local ring then $D$ is too. If $K$ is a field of characteristic   zero and $D$ is a $K$-algebra, regular local ring, essentially of finite type over $K$ then $D$ is smooth over $K$. 

The Zariski's Uniformization Theorem could have the following form

\begin{Theorem} \label{T1}  (Zariski \cite{Z}) A valuation ring containing a field $K$ of characteristic zero is a    a filtered union of its smooth $K$-subalgebras.
\end{Theorem}

   Let $k$ be a field of characteristic $2$ and $K=k(X)$ the field of fractions of the polynomial ring $k[X]$. Then $E=K[Y]/(Y^2-X)$ is a field, so a regular local ring,
   but $E/K$ is not separable, so $ E $ is not a smooth $K$-algebra. Thus in positive characteristic a similar statement of the above theorem needs to say more then the fact that a valuation ring containing a field $K$ is a filtered union of its regular local $K$-subrings.
   
 However, B. Antieau, R. Datta \cite[Theorem 4.1.1]{AD} showed the following theorem. 
\begin{Theorem} (\cite{AD})\label{T2}
  Every perfect valuation ring of characteristic $p>0$ is a filtered union of its smooth ${\bf F}_p$-subalgebras. 
 \end{Theorem}
 
 Recently, we showed the following theorem.
 
   \begin{Theorem}(\cite[Corollary 5]{P3})\label{T3} Let $V$ be a  valuation ring containing its residue field $k$   with a $\bf Z$-free value group $\Gamma$ (for example when it is finitely generated) and $K$ its fraction field. Assume that  $K=k(x,y)$ for some  algebraically independent elements $x=(x_i)_{i\in I},y=(y_j)_{j\in J}$ over $k$ such that $\val(y)$ is a basis in $\Gamma$. 
 Then $V$ is a filtered  union of its smooth $k$-subalgebras. In particular, $V$ is a filtered  union of its  $k$-subalgebras, which are regular local rings, essentially of finite type over $k$.
\end{Theorem}

It is the purpose of this paper to explain this result together with some of its applications.

\section{Immediate extension of valuation rings}

An {\it{immediate extension}} of valuation rings is an extension inducing trivial extensions on residue fields and group value extensions.

Using \cite[Theorem 1,VI,(10.3)]{Bou} we got in \cite[Lemma 26 (1)]{P} the following result.
   \begin{Lemma}(\cite{P})\label{l} Let $V$ be a  valuation ring containing its residue field $k$    with a $\bf Z$-free value group $\Gamma$ (for example when it is finitely generated), $\val$ its valuation and $K$ its fraction field. Assume that  $K=k(x)$ for some   elements $x=(x_i)_{i\in I}$ over $k$ such that $\val(x)$ is a basis in $\Gamma$. 
 Then $V$ is a filtered  union of its  $k$-subalgebras, which are localizations of some polynomial rings over $k$ in some $x$, so some smooth $k$-algebras.
\end{Lemma}

Let $V$ be a valuation ring  with a $\bf Z$-free  value group  $\Gamma$,  containing its residue field $k$  and $x$ some elements of $V$ such that $\val(x)$ is a basis of $\Gamma$.
Our idea was to split in two parts  a proof of a possible extension of Zariski's Uniformization Theorem in positive characteristic. Firstly we noticed that $W=V\cap k(x)$ is a  a filtered  union of its smooth $k$-subalgebras by Lemma \ref{l}. Then we try to find cases when the immediate extension $V/W$ is  a  a filtered  union of its smooth $W$-subalgebras. Unfortunately, in general  this is not the case (see \cite[Example 3.1.3]{Po1},\cite{O},  \cite[Example 3.6]{P4}) and for some positive results we have to state different other necessary results.

The above idea was used already in \cite{P}, where we proved the following result.

\begin{Theorem} (\cite[Theorem 2]{P})
 If $V\subset V'$ is an immediate extension of valuation rings   containing $\bf Q$
then $V'$ is a filtered direct colimit of smooth $V$-algebras. 
\end{Theorem}

Let $V$ be a valuation ring, $\lambda$ be a fixed limit ordinal  and $v=\{v_i \}_{i < \lambda}$ a sequence of elements in $V$ indexed by the ordinals $i$ less than  $\lambda$. Then $v$ is called \emph{pseudo convergent} if 

$\val(v_{i} - v_{i''} ) < \val(v_{i'} - v_{i''} )     \ \ \mbox{for} \ \ i < i' < i'' < \lambda$
(see \cite{Kap}, \cite{Sch}).
A  \emph{pseudo limit} of $v$  is an element $w \in V$ with 

$ \val(w - v_{i}) < \val(w - v_{i'}) \ \ \mbox{(that is,} \ \ \val(w -  v_{i}) = \val(v_{i} - v_{i'}) \ \ \mbox{for} \ \ i < i' < \lambda$.

k
\begin{Lemma}(\cite[Proposition 16]{P3}) \label{l1}
Let $V \subset V'$ be an  immediate extension of valuation rings,  $K, K'$ the fraction fields of $V, V'$, $\mm$ the maximal ideal of $V$ and  $y\in K'$ an element  which is not in  $K$. Assume that  $y$ is  a pseudo limit of a pseudo convergent sequence $v=(v_j)_{1\leq j<\lambda}$ over $V$, which has no pseudo limit in $K$.
 Set $y_j=(y-v_j)/(v_{j+1}-v_j)$.
Then for every nonzero polynomial  $g =\sum_{e=0}^m a_e Y^e\in V[Y]$ such that  $a_e=0$ for any $e>o$ multiple of $p$ if char\ $K=p>0$,  and every ordinal $1\leq \nu<\lambda$,  there exist some  $\nu<t<\lambda$ and a polynomial $g_1\in V[Y_t]$ such that 
$g(y)=g_1(y_t)$
 and    $g_1=g_1(0)+c  Y_t+g_2,$
for some $c\in V\setminus \{0\}$ and  $g_2\in  c\mm Y_t^2 V[Y_t]$. 
\end{Lemma}

For the proof very important is that changing $Y\to (Y-v_j)/(v_{j+1}-v_j)$ in $g$  we may assume that the nonzero coefficients of $g-g(0)$ have different values for $j$ large enough. It is worth to note that $\val(g(y)) $ could be different from $\val(g(y_j))$ for all $j$ as shows \cite[Example 16]{P4}. 

\begin{Remark} \label{r}
If the characteristic of $V$ is $p>0$ then taking above $g\in V[Y^p]$ we see that for any $t$ the polynomial $g_1=g(Y_t)\in V[Y_t^p]$ and so the coefficient of $Y_t$ in $g_1$ is zero.
\end{Remark}

When $p>0$ and there exist  in $g$ above some nonzero coefficients of some $Y^n$ with $n>0$ multiple of $p$ then we get the following lemma.

\begin{Lemma}(\cite[Corollary 17]{P3})\label{l2}
Let $V \subset V'$ be an  immediate extension of valuation rings,  $K, K'$ the fraction fields of $V, V'$, $\mm$ the maximal ideal of $V$ and  $y\in K'$ an  element  which is not in  $K$. 
Assume that char\ $K=p>0$ and $y$ is  a pseudo limit of a pseudo convergent sequence $v=(v_j)_{1\leq j<\lambda}$ over $V$, which  has no pseudo limit in $K$.
 Set $y_j=(y-v_j)/(v_{j+1}-v_j)$.
Then for every nonzero  polynomial  $g\in V[Y]$  and every ordinal $1\leq \nu<\lambda$ one of the following statements holds.

\begin{enumerate}
\item  There exist some  $\nu<t<\lambda$ and a polynomial $g_1\in V[Y_t]$ such that 
$g(y)=g_1(y_t)$
 and    $g_1=g_1(0)+c  Y_t+g_2,$ 
for some $c\in V\setminus \{0\}$ and  $g_2\in  c\mm Y_t^2 V[Y_t]$.
\item 
There exist some  $\nu<t<\lambda$ and a polynomial $g_1\in V[Y_t]$ such that 
$y g(y)=g_1(y_t)$
 and    $g_1=g_1(0)+c  Y_t+g_2,$ 
for some $c\in V\setminus \{0\}$ and  $g_2\in  c\mm Y_t^2 V[Y_t]$.
\end{enumerate}
\end{Lemma}

These lemmas are applied  to Lemma \ref{l3} given below together with the next lemma which is in fact our key lemma.

\begin{Lemma}(\cite[Lemma 21]{P3}) \label{k}
Let $V \subset V'$ be an  immediate extension of valuation rings,  $K, K'$ the fraction fields of $V, V'$, $\mm, \mm'$ the maximal ideals of $V,V'$, 
$(y_e)_{1\leq e\leq  m}$ some elements of $V'$ and $g_e\in V[Y_1,\ldots,Y_m]$, $1\leq e<m$ some polynomials such that the determinant of $((\partial g_e/\partial Y_i)((y_{e'}))_{1\leq e< m, 1\leq i<m}$ is not in $\mm'$.
 Assume that $y_m$ is transcendental over $K$ and $g_e((y_{e'})_{e'})=0$ for all $1\leq e<m$.
 Then $V[y_1,\ldots,y_m]_{\mm'\cap V[y_1,\ldots,y_m]}$ is  a smooth $V$-subalgebra  of $V'$.
\end{Lemma}

\begin{Lemma}(\cite[Lemma 22]{P3}) \label{l3}
 Let  $ V'$ be an  immediate extension  of a valuation ring $V$ containing a field, $\mm, \mm'$ the maximal ideals of $V,V'$, $K\subset K'$ their fraction field extension and $y\in V'$ an unit, transcendental element over $K$. 
Let  $f\in V[Y]$ be a nonzero polynomial and some $d\in V\setminus \{0\}$ and an unit $z\in V'$ such that $f(y)=d z$. Then there exists 
 a smooth $V$-subalgebra  of $V'$ containing $y,z$.  
\end{Lemma}

The above lemma gives an  idea of the proof of our  desingularization (in fact smoothification).
By \cite[Theorem 1]{Kap}  $y$ is a pseudo limit of a pseudo convergent sequence $v=(v_j)_{j<\lambda}$, which 
 has no pseudo limits in $K$.  Set $y_j=(y-v_j)/(v_{j+1}-v_j)$.

 By Lemmas \ref{l1}, \ref{l2} there exist some  $j<\lambda$ and a polynomial $h\in V[Y_j]$ such that either 
$f(y)=h(y_j)$, or $y f(y)=h(y_j)$
 and  $h=h(0)+c  Y_j+h',$
for some   $c\in V\setminus \{0\}$ and $h'\in  c\mm Y_j^2 V[Y_j]$.   If $c|d$ then $c|h(0)$ and $(h-dZ)/c$, or $(h-dYZ)/c$ is in the kernel of the map $V[Y_j,Z]\to V'$,
$(Y_j,Z)\to (y_j,z)$ and  $V[y_j,z]_{\mm'\cap V[y_j,z]}$ is smooth over $V$ by Lemma \ref{k} containing  $y,z$.

\section{A smoothification theorem}

Mainly we try to explain here the following theorem

\begin{Theorem}(\cite[Theorem 23]{P3}) \label{T4}
 Let  $ V'$ be an  immediate extension  of a valuation ring $V$ containing a field, $\mm, \mm'$ the maximal ideals of $V,V'$, $K\subset K'$ their fraction field extension and $y_0\in V'$ an unit, transcendental element over $K$ with $K'=K(y_0)$. 
  Then $V'$ is  a filtered  union of its  smooth  $V$-subalgebras.
\end{Theorem}

For the proof we need also some extensions of Lemmas \ref{l1}, \ref{l2} given for two variables. Their proofs are similar but  more complicated.

\begin{Lemma}(\cite[Proposition 18]{P3}) \label{l4}
Let $V \subset V'$ be an  immediate extension of valuation rings,  $K, K'$ the fraction fields of $V, V'$, $\mm$ the maximal ideal of $V$ and  $y_1,y_2\in K'$ two elements,  which are not in $K$. Assume that $y_i$, $i=1,2$ is  a pseudo limit of a pseudo convergent sequence $v_i=(v_{i,j})_{1\leq j<\lambda_i}$, $i=1,2$ over $V$, which   have no pseudo limit in $K$.
 Set $y_{ij}=(y_i-v_{i,j})/(v_{i,j+1}-v_{i,j})$.
Then for every nonzero polynomial   $g=\sum_{e_1,e_2=0}^m a_{e_1,e_2} Y_1^{e_1} Y_2^{e_2} \in V[Y_1,Y_2]$, such that $a_{e_1,e_2}=0$ if char\ $K=p>0$ and at least one of $e_1,e_2$ is a nonzero multiple of  $p$,  and every ordinals $1\leq \nu_i<\lambda_i$,  $i=1,2$  there exist some  $\nu_i<t_i<\lambda_i$, $i=1,2$ and a polynomial $g_1\in V[Y_{1,t_1},Y_{2,t_2}]$ such that 
$g(y_1,y_2)=g_1(y_{1,t_1},y_{2,t_2})$
 and  $g_1=g_1(0)+c_1  Y_{1,t_1}+c_2 Y_{2,t_2}+g_2,$
for some   $c_1,c_2\in V, $ at least one of them nonzero, and $g_2\in  (c_1,c_2)\mm (Y_{1,t_1},Y_{2,t_2})^2 V[Y_{1,t_1},Y_{2,t_2}]$.   
\end{Lemma}

\begin{Lemma}(\cite[Corollary 20]{P3}) \label{l5}
Let $V \subset V'$ be an  immediate extension of valuation rings,  $K, K'$ the fraction fields of $V, V'$, $\mm$ the maximal ideal of $V$ and  $y_1,y_2\in K'$ two  elements  which are not in  $K$. 
Assume that char\ $K=p>0$ and $y_i$, $i=1,2$ are   pseudo limits of two pseudo convergent sequences $v_i=(v_{i,j})_{1\leq j<\lambda_i}$, $i=1,2$ over $V$, which have no pseudo limits in $K$.
 Set $y_{i,j}=(y_i-v_{i,j})/(v_{i,j+1}-v_{i,j})$.
Then for every nonzero polynomial  $f\in V[Y_1,Y_2]$  and every two ordinals $\nu_i<\lambda_i$ 
there exist some  $1\leq \nu_i<t_i<\lambda_i$, $i=1,2$ and a polynomial $g\in V[Y_{1,t_1}, Y_{2,t_2}]$
 such that $g=g(0)+c_1  Y_{1,t_1}+c_2  Y_{2,t_2}+g',$ 
for some $c_1,c_2\in V$, at least one of them nonzero,   $g'\in  (c_1,c_2)\mm (Y_{1,t_1},Y_{2,t_2})^2   V[Y_{1,t_1}, Y_{2,t_2}]$ 
and $g(y_{1,t_1},y_{2,t_2})$ is one of the following elements
$f(y_1,y_2)$, \ $y_1f(y_1,y_2)$, \ $y_2f(y_1,y_2)$, \ $y_1y_2f(y_1,y_2)$.
\end{Lemma}
The proof of Theorem \ref{T4} goes by induction  using Lemmas \ref{l1}, \ref{l2},   \ref{k}, \ref{l4}, \ref{l5}
as  in the proof of Lemma \ref{l3}. If
 $f_i(y)=d_iz_i$, $i=1,\ldots, s$
for some $f_i\in V[Y]$, $d_i\in V\setminus \{0\}$, $z_i\in V'$ then there exists a smooth $V$-subalgebra of $V'$  containing $y,(z_i)$.

Using by induction Theorem \ref{T4} we can extend it  for the case when $V\subset V'$ is a pure transcendental extension.
\begin{Theorem}(\cite[Theorem 4]{P3})\label{T5} Let  $V\subset V'$ be an immediate  extension of valuation rings   containing  a field and $K\subset K'$ its fraction field extension. If $K'=K(x)$ for some algebraically independent system of elements $x$ over $K$ 
then $V'$ is a filtered union of smooth $V$-subalgebras of $V'$. 
\end{Theorem}

{\bf Proof of Theorem \ref{T3}} 

Use  Theorem \ref{T7} and Lemma \ref{l}.

\section{Some applications}

A {\em Henselian} local ring is a local ring 
 $(A,\mm)$ for which the Implicit Function Theorem hods, that is 
 for every system of polynomials  $f=(f_1,\ldots ,f_r) $ of $A[Y]$,
$Y=(Y_1,\ldots,Y_N)$, $N\geq r$ over  $A$ and every solution  $y=(y_1,\ldots,y_N)\in A^N$     modulo $\mm$ of $f$ in   $A$ such that an $r\times r$-minor $M$ of Jacobian matrix $(\partial f_i/\partial Y_j)_{1\leq i\leq r,1\leq j\leq N}$  satisfies $M(y)\not \in \mm$ there exists a solution $\tilde y$ of $f=0$ in  $A$ with $y\equiv {\tilde y}$ modulo $\mm$.
 Thus  $(A,\mm)$ is Henselian if for any smooth $A$-algebra $B$  any $A$-morphism $B\to A/\mm$ can be lifted to an $A$-morphism $w:B\to A$. Thus  $w$ is a retraction of $A\to B$.

As in \cite[Proposition 18]{P1}(see also \cite[Corollary 24]{P3}, \cite[Proposition 1]{P4}) we have the following corollary.

\begin{Corollary}(\cite{P4})\label{C} Let $V\subset V'$ be an immediate extension of valuation rings containing a field and $K\subset K'$ its fraction field extension. If $K'=K(x)$ for some algebraically independent system of elements $x$ over $K$ 
and $V$ is Henselian then every system of polynomials over $V$, which has a solution in $V' $ has also one in $V$.
\end{Corollary}

The ideas of this corollary were used to show some conjectures of M. Artin from \cite{Ar} in 
\cite[Theorems 1.3, 1.4]{Po0} (see also \cite[Theorem 5.3.1]{I}) using the so called the General N\'eron desingularization (see  \cite[Theorem 2.5]{Po0}, \cite[Theorem 1.1]{S} and \cite[Theorem 5.2.56]{I}).

For a general immediate extension of valuation rings $V\subset V'$ we should expect that $V'$ is a  filtered  union of its complete intersection $V$-subalgebras. 
 
A {\em complete intersection} $V$-algebra {\it{essentially of finite type}} is a local $V$-algebra of type $C/(P)$, where $C$ is a localization of a polynomial $V$-algebra of finite type and $P$ is a regular sequence of elements of $C$. 

\begin{Theorem}  (\cite[Theorem 1]{P2})\label{T6} Let  $ V'$ be an  immediate extension  of a valuation ring $V$  and  $K\subset K'$ the fraction field extension. If $K'/K$ is algebraic  then $V'$ is a filtered
 union of its complete intersection $V$-subalgebras of finite type.
\end{Theorem}

As a consequence of Theorems \ref{T5}, \ref{T6} the next theorem follows.

\begin{Theorem} (\cite[Theorem 6]{P3}) \label{T7} Let  $ V'$ be an  immediate extension  of a valuation ring $V$  and  $K\subset K'$ the fraction field extension. Then $V'$ is a filtered
 union of its complete intersection $V$-subalgebras of finite type.
\end{Theorem}

\begin{Corollary} (\cite[Corollary 7] {P3})Let $V$ be a  valuation ring containing its residue field $k$   with a $\bf Z$-free value group $\Gamma$. 
 Then $V$ is a filtered  union of its complete intersection $k$-subalgebras. 
\end{Corollary}
For the proof choose  some elements $y$ in $V$ such that $\val(y)$ is a basis in $\Gamma$  and apply Lemma \ref{l}   for $V\cap k(y) $ and the above theorem for the immediate extension $V\cap k(y)\subset V$.

 \end{document}